\documentclass[12pt]{article}

\usepackage{amssymb}
\usepackage{amsmath}
\usepackage{graphicx}
\usepackage{indentfirst}

\begin{document}

\newtheorem{Theorem}{Theorem}[section]
\newtheorem{thm}[Theorem]{Theorem}
  \newtheorem{cor}[Theorem]{Corollary}
 \newtheorem{lem}[Theorem]{Lemma}
 \newtheorem{prop}[Theorem]{Proposition}
  \newtheorem{defn}[Theorem]{Definition}
  \newtheorem{rem}[Theorem]{Remark}
  \newtheorem{algm}[Theorem]{Algorithm}
  \newtheorem{conj}[Theorem]{Conjecture}
\newtheorem{exa}[Theorem]{Example}
\newtheorem{prob}[Theorem]{Problem}
 \def\dfrac#1#2{{\displaystyle{#1\over#2}}}

\def\I{{\bf i}}
\def\G{{G^{\sigma}}}
\def\Gr{\overrightarrow{G}}
\def\Hr{\overrightarrow{H}}
\def\H{{H^{\tau}}}
\def\T{\overrightarrow{T}}
\def\Gl{\overleftarrow{G}}
\def\P{\overrightarrow{P}}
\def\a{\lambda}
\def\zm{\noindent{\bf Proof.  }}
\def\ezm{ $ \hfill \Box$ \vspace*{0.2cm}}
\def\Z{\mathbb{Z}}

\thispagestyle{empty}

\title{The graphs  with  all but  two eigenvalues equal to $2$ or $-1$ }
\author{Jing Li, Deqiong Li and Yaoping Hou\thanks{Corresponding author: yphou@hunnu.edu.cn}  \\[2mm]
{\small Department  of Mathematics,  Hunan    Normal University} \\
{\small Changsha, Hunan 410081, China} }

\maketitle
\vspace*{0.15cm}
\begin{abstract}
In this paper,  all graphs whose adjacency matrix has at most two eigenvalues (multiplicities included) different from $2$ and $-1$  are determined. These graphs
conclude a class of generalized friendship graphs $F_{t,r,k}, $ which is the graph of $k$ copies of the complete graph $K_t$  meeting in common $r$ vertices such that $t-r=3.$ Which of these graphs are determined by its  spectrum  is are also obtained.
 \end{abstract}


\normalsize
\baselineskip=0.25in
\section{ Introduction}

All graphs in this paper are simple graphs and all spectrum of a graph are adjacency spectrum. Let $G=(V,E)$ be a graph. The adjacency matrix $A(G)$ (or $A$) of $G$ is an $n\times n$ matrix, whose $(i,j)$-entry is $1$ if vertex $v_{i}$ is adjacent to $v_{j}$ (denote by $v_i \sim v_j$),  and is $0$ otherwise. The characteristic polynomial $P_G(x)=det(xI-A(G))$  is called the characteristic polynomial of $G$. The eigenvalues of $A$  are called the
adjacency eigenvalues  of $G$. There are many results on the eigenvalues of graphs and their application, see \cite{b1} for more details.

 Connected graphs with a small number of distinct eigenvalues have aroused a lot of interest in the past several decades. This problem was first raised by Doob \cite{b13}. It is well known that a connected graph with just two distinct  eigenvalues if and only if it is completed graph and a regular connected graph with just three distinct eigenvalues if and only if it is strongly regular graph. It is difficult to characterise all non-regular connected graphs with three or four distinct eigenvalues. There are interesting results on regular graphs with four distinct eigenvalues \cite{b10},  non-regular graph with three distinct eigenvalues \cite{b11}, biregular graphs with three distinct eigenvalues \cite{b9} and small regular with four distinct eigenvalues \cite{b12}.  Cioab\v{a} et al. in \cite{b8} determined all connected graphs with at most two eigenvalues  different from $-2$ or $0$.
  For more results on graphs with few  distinct eigenvalues, we refer the reader to \cite{b15,b14,b16}.

For $0\leq r \leq t$, denote the generalized friendship graph on $kt-tr+r$ vertices by $F_{t,r,k},$ where $F_{t,r,k}$  is the graph of $k$ copies of  the complete graph $K_t$  meeting in a common $r$ vertices. Clearly $F_{t,r,1}=F_{t,t,k}=K_t,$ which is determined by its spectrum. For convenience we shall assume that $k\geq 2.$  $F_{3,1,k}$ is the  friendship graph, which  is determined by its spectrum if $k\not=16$  \cite{b5}. It is not difficult to obtain that the spectrum of  $F_{t,r,k}$ has at most two   eigenvalues (multiplicities  included) different from $t-r-1$ and $-1.$  It may be a interesting problem  that   $F_{t,r,k}$  is  whether determined by its  spectrum.
Very recently, Cioab\v{a} et al. in \cite{b5} determined all connected graphs with at most two eigenvalues  different from $\pm1,$  which responds to the case $t-r=2,$ and prove that
friendship graph $F_{3,1,k}$ is  determined by its spectrum unless $k=16.$

In this paper, we consider the case of $t-r=3$  and determine all connected graphs  with  two eigenvalues  different from $2$ and $-1$, these graphs consist of four infinite families and twenty sporadic graphs, which of these graphs are determined by its  spectrum is  also obtained.

\section{Main tools}
\par We start with a well known result on equitable partitions (see for example \cite{b1} ). Consider a partition $\mathcal{P} = \{V_1,\dots , V_m\}$ of the set $V= \{1, \dots , n\}$. The characteristic matrix $\mathcal{X}_\mathcal{P}$ of $\mathcal{P}$ is the $n \times m$ matrix whose columns are the character vectors of $V_1, \dots, V_m$. Consider a symmetric matrix $A$ of order $n$, with rows and columns partitioned according to $\mathcal{P}$. The partition of $A$ is equitable if each submatrix  $A_{i, j}$  formed by the rows of $V_{i}$ and the columns of $V_{j}$ has constant row sums $q_{ij}$. The $m\times m$ matrix $Q =(q_{i,j})$ is called the quotient matrix of $A$ with respect to $\mathcal{P}$.
\begin{lem}\label{odd}\cite{b1}
The matrix $A$ has the following two kinds of eigenvectors and eigenvalues:

(1) The eigenvectors in the column space of $\mathcal{X}_\mathcal{P}$; the corresponding eigenvalues coincide with the eigenvalues of $Q$;

(2) The eigenvectors orthogonal to the columns of $\mathcal{X}_\mathcal{P}$; the corresponding eigenvalues of $A$ remain unchanged if some scalar multiple of the all-one block $J$ is added to block $A_{i,j}$ for each $i,j\in \{{1,\dots, m}\}$.\end{lem}
\par   The degree of a vertex $v$, denoted by $d_v$, which is the number of vertices adjacent to $v$, $d_{uv}$ is the number of common neighbors of $u$ and $v$.  If  the vertices $i$ and $j$ are adjacent, we denoted by $i\sim j$, otherwise  $i\nsim j$. Let $mK_{3}$ denote the disjoint union of $m$ triangles, and $kK_{2}$ denote the disjoint union of $k$ edges, and $T_{3m}$ be the adjacency matrix of $ mK_{3}$ and $R_{2k}$ be the adjacency matrix of $ k K_2$. We denote the $m\times n$ all-ones matrix by $J_{m,n}$ (or just $J$ ) and the  $m\times n$ all-zeros matrix by $0_{m}$ (or $0$). We define a  $2k \times k$ matrix  $S_{2k}$ as following:
$$ S_{2k}=\begin{bmatrix}\begin{smallmatrix}
                          1&0&0&\cdots&0&0\\
                          1&0&0&\cdots&0&0\\
                          0&1&0&\cdots&0&0\\
                          0&1&0&\cdots&0&0\\
                          \vdots&\vdots&\vdots& &\vdots&\vdots\\
                          0&0&0&\cdots&0&1\\
                          0&0&0&\cdots&0&1\\
  \end{smallmatrix}
\end{bmatrix}.
$$
\begin{lem}\label{even}\cite{b1}
Let $G$ be a graph with smallest eigenvalue $-1$, then $G$ is the disjoint union of complete graphs.
 \end{lem}

 \begin{lem}\label{clique}(\cite{b6})
 The only connected graphs having the largest eigenvalue $2$ are the graphs in Figure 1 .\end{lem}
 \begin{figure}[ht]
  \centering
  \includegraphics[width=10cm, height=3cm]{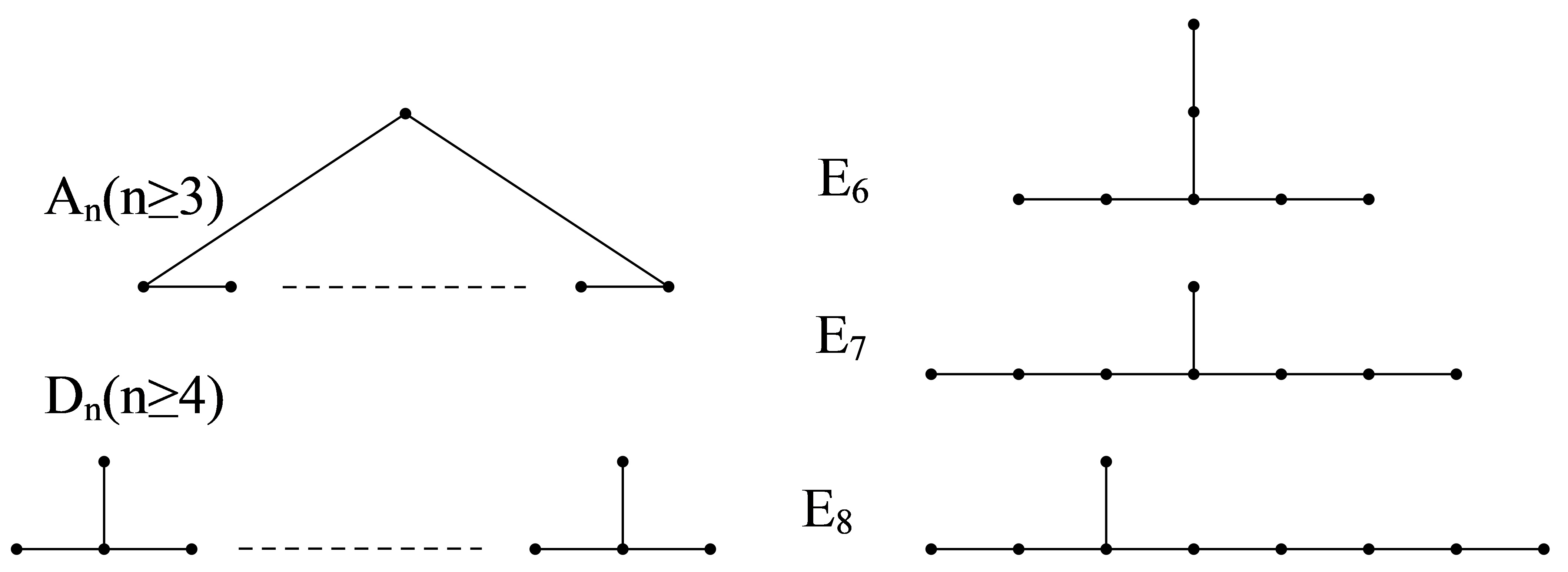}\\
  \caption{Connected graphs with the largest eigenvalue 2. }\label{fig-1}
\end{figure}

  \begin{prop} \label{p:wAPPROX}
Let $G$ be a graph with $n$ vertices, we have

(i) If $G$ has all its eigenvalues equal to $2$ and $-1$, then $G = \frac{n}{3} K_3$.

(ii) If $G$ has all but one eigenvalue equal to $2$ and $-1$, then $G$ is the disjoint union of complete
graphs with all but one connected components equal to $K_3$.

(iii) If $G$ has just two eigenvalues, $r$ and $s$ ($r \geq s$) different from $2$ and $-1$, then $r > 2$ and
$s < -1$, or $G$ is a disjoint union of complete graphs with two connected components different from $K_3$.
\end{prop}
\zm If $G$ has the smallest  eigenvalue $-1$, by Lemma \ref{even}, then $G$ is the disjoint union of complete graphs, which leads to (i),(ii) and the second option of (iii). If $G$ has the largest eigenvalue $ 2 $, by Lemma \ref{clique}, then $G$ are the graphs in Figure 1. Computing eigenvalue of these graphs, the corresponding graphs are not in $G$, therefore $r>2$, and $s<-1$, this case is captured by the first option of (iii).
\ezm

 By Proposition  \ref{p:wAPPROX},  in order to obtain the connected graphs with at most two eigenvalues differen from 2 and $-1,$ it is sufficient to  determine the graphs with  just two eigenvalues $r$ and $s$ ($r >2>-1 > s$) different from $2$ and $-1.$ Therefore, the spectrum of such a graph $G$ has two interesting properties: The first property is that the second largest eigenvalue of $A(G)$ is $2$, and the second smallest eigenvalue is equal to $-1$. By eigenvalue interlacing, this gives a considerable reduction on the possible induced subgraphs (see Lemma \ref{main}). The second property is that $(A(G)+I)(A(G)-2I)$ has rank $2$ and is positive semi-definite. This leads to conditions for the structure of $(A(G)+I)(A(G)-2I)$ (see Lemmas \ref{nonsigular}, \ref{nine}). Because of these observations, we take a more general approach,  and consider all graphs with the mentioned two properties.  In what followings we determine all connected graphs with only two eigenvalues $r$ and $s$ $(r>2>-1>s)$ different from $2$ and $-1$.

 \begin{lem}\label{nonsigular}
If the graph $G$ with only two eigenvalues $r>2$ and $s<-1$ (multiplicities included) different from $2$ and $-1$, then

(i) One connected component of $G$ has all vertices with  degree at least $3$, and all other connected components are isomorphic to $K_3$.

(ii) If the vertices $u\nsim v$, and each neighbor of $u$ is also a neighbor of $v$, then $d_v - d_u \geq 5$.
 \end{lem}

 \zm (i) We prove the result by contradiction, suppose $u$ is a vertex of degree $1$, $v$ is a vertex of degree $2$. Let $v$ be the neighbor of $u$, and assume that $v$ has another neighbor $w$ of degree $d_w$. The $2\times2$ principal submatrix of $A^2-A-2I$ corresponding to $u$ and $w$ equals
 $$S=\left[
                        \begin{array}{cc}
                          -1 & 1  \\
                          1  & d_w-2 \\
                        \end{array}
                      \right].
 $$
  The $2\times2$ principal submatrix of $A^2-A-2I$ corresponding to $v$ and $w$ equals $$S'=\left[
                        \begin{array}{cc}
                          0 & -1  \\
                          -1 & d_w-2 \\
                        \end{array}
                      \right].
 $$

  We have $\det S<0$, det $S'<0$, which contradicts with that $A^2-A-2I$ is positive semi-definite. Thus we have $d_x\geq 3$ for any vertex $x\in G$.

 (ii) The $2\times2$  principal submatrix of $A^2-A-2I$ corresponding to $u$ and $v$ equals
 $$S=\left[
                        \begin{array}{cccc}
                          d_u-2 & d_u  \\
                          d_u & d_v-2 \\
                        \end{array}
                      \right].
 $$
If $d_v\leq d_u + 4$, then det $ S\leq (d_u-2)(d_u + 2)-d_u^2<0$, contradiction. \ezm
\par Note that  Lemma \ref{nonsigular} (ii) indicates that any two non-adjacent vertices can not have the same set of neighbors.

\begin{lem}\label{five}\cite{b1} Let $G$ be a bipartite graph, if $\lambda$ is an eigenvalue of $G$ with multiplicity $k$, then $-\lambda$ is also an eigenvalue of $G$ with multiplicity $k$.
\end{lem}

\begin{lem}\label{prod} (Interlacing Theorem)\cite{b1}
Let $A$ be a symmetric $n\times n$ matrix and let $B$ be a principal submatrix of $A$ of order $n-1$. If $\lambda_1\geq\dots\geq\lambda_n$ and $\mu_1\geq\dots\geq\mu_{n-1}$ are the eigenvalues $A$ and $B$, respectively, then
$$\lambda_1\geq\mu_1\geq\lambda_2\geq\dots\geq\lambda_{n-1}\geq\mu_{n-1}\geq\lambda_n.$$
\end{lem}

\begin{figure}[h]
  \centering
  \includegraphics[width=13cm, height=12cm]{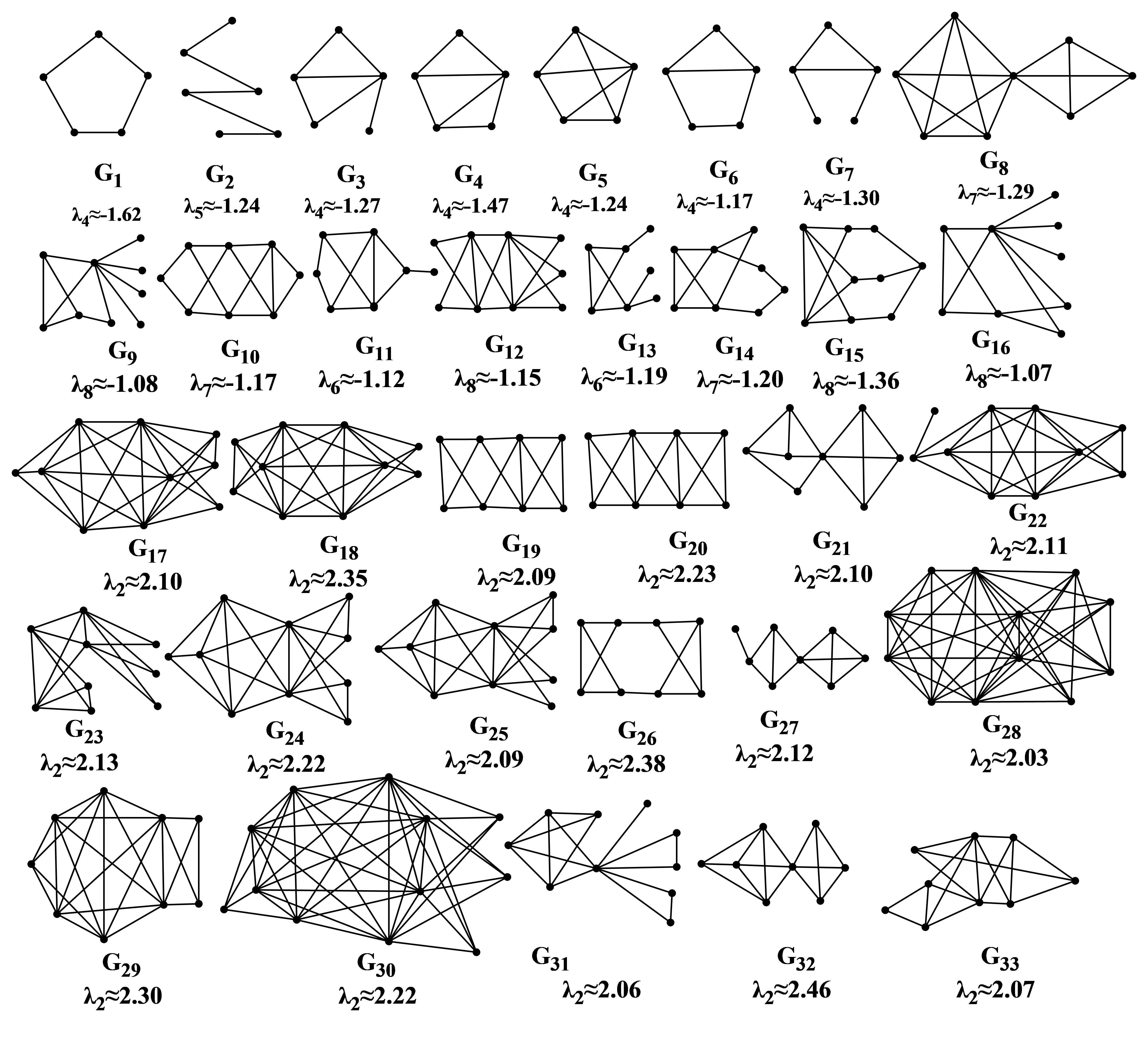}\\
  \caption{Forbidden induced subgraphs.}\label{fig_2}
\end{figure}

 Define $\mathcal{F}$ to be the set of connected graphs with two eigenvalues $r > 2$ and $s < -1$ (multiplicities included), and all other eigenvalues equal to $2$ and $-1$. Lemmas \ref{even}, \ref{five} indicate that the graph $G\in \mathcal{F}$ is not bipartite.  In order to find all graphs with only two eigenvalues different from $2$ and $-1$, we start with a list of forbidden induced subgraphs.

 \begin{lem}\label{main}   No graph in $\mathcal{F}$ has one of the graphs presented in Figure $2$ as an induced subgraph.
 \end{lem}
 \zm Each graph in Figure $2$  has its second largest eigenvalue $\lambda_2$ strictly greater than $2$, or its second smallest eigenvalue $\lambda_{n-1}$ strictly less than $-1$. Interlacing completes the proof.
 \ezm

\section{Main results}
We begin with the description of the graphs in $\mathcal{F} $. The proof will be given in the next section.
\begin{thm}\label{one} For each $G\in \mathcal{F} $, the adjacency matrices and the corresponding spectra of $G$ are one of the following forms:
 \vskip 0.10 cm
 (i). $\begin{bmatrix}
                          J-I_a & J  \\
                          J & T_{3k} \\
                     \end{bmatrix} (a\geq 1,k\geq2)$
                       with spectrum $\{{\frac{(a+1)\pm\sqrt{(a-3)^2+12ak}}{2}, 2^{k-1}, -1^{2k+a-1}}\},$ \vskip 0.20 cm

(ii). $\begin{bmatrix}
                          T_{3k} & J   \\
                          J & T_{3\ell} \\
                      \end{bmatrix}(k\geq \ell\geq2)$
 with spectrum $\{{2\pm3\sqrt{k\ell}}, 2^{k+\ell-2}, -1^{2(k+\ell)}\},$    \vskip 0.20 cm

(iii). $\left[
      \begin{array}{cc}
       R_{2m }& J-S_{2m}  \\
       J-S_{2m}^{T} & 0 \\
      \end{array}
    \right] (m\geq3)$  with spectrum $\{{\frac{1\pm\sqrt{9-16m+8m^{2}}}{2}, 2^{m-1}, -1^{2m-1}}\},$   \vskip 0.20 cm

(iv).
$\begin{bmatrix}
                          J-I_6 & J  & 0  \\
                          J & T_{3k}  & J \\
                          0 & J & R_2  \\
\end{bmatrix} (k\geq2)
 $ with spectrum $\{{{3\pm2\sqrt{1+6k}}, 2^{k}, -1^{2k+6}}\},$ \vskip 0.20 cm
(v). $\begin{bmatrix}
                          J-I_a & J  & J \\
                          J & J-I_b  & 0 \\
                          J & 0 & R_{2} \\
                        \end{bmatrix}$ where $(a,b)=(2,9), (3,6)$ and $(6,5)$,
  with the corresponding  spectra $\{{4\pm\sqrt{37}},2,-1^{10}\}$,$\{3\pm2\sqrt{7},2,-1^{8}\}$, $\{4\pm3\sqrt{5},2,-1^{10}\},$  \vskip 0.20 cm

(vi). $\begin{bmatrix}
                          J-I_a & J  & J \\
                          J & J-I_b  & 0 \\
                          J & 0 & 0 \\
                     \end{bmatrix}
                        $
  where $(a,b)=(7,45),(8,27),(9,21),(10,18),
  (12,15),$\vskip 0.25 cm $(15,13)$, $(18,12),(24,11)$ and $(42,10),$ with the corresponding spectra $$\{{24\pm\sqrt{730}},2,-1^{50}\},\{\frac{31\pm9\sqrt{17}}{2},2,-1^{33}\},\{13\pm\sqrt{259},2,-1^{28}\},$$
     $$\{{12\pm\sqrt{229}},2,-1^{26}\},\{\frac{23\pm\sqrt{865}}{2},2,-1^{25}\},\{12\pm3\sqrt{26},2,-1^{26}\},$$
     $$\{{13\pm2\sqrt{67}},2,-1^{28}\},  \{\frac{31\pm\sqrt{1441}}{2},2,-1^{33}\},\{24\pm3\sqrt{85},2,-1^{50}\},$$

(vii). $\begin{bmatrix}
J-I_a & J& 0 \\
  J & 0   & J-S^{T}_{2m}   \\
  0 & J-S_{2m} & R_{2m} \\
\end{bmatrix}$
where $(a,m)=(4,4)$ and $(6,3),$ \\ with corresponding spectra $\{7,-5,2^4,-1^{10}\}$  and  $\{2\pm \sqrt{33},2^3,-1^{10}\}$.

(viii). $\begin{bmatrix}
                          J-I_a & J & 0 \\
                          J & R_{2k} &  J-S_{2k}  \\
                          0 & J-S^{T}_{2k} & 0 \\
                       \end{bmatrix}$
                      where $(a,k)=(4,10),(5,7), (6,6)$ and $(9,5),$ \\ with the corresponding spectra
   $\{1\pm2\sqrt{61},2^{10},-1^{22}\}$, $\{\frac{3\pm3\sqrt{65}}{2},2^{7},-1^{17}\},$\\ $\{2\pm\sqrt{129},2^{6},-1^{16}\}$ and  $\{\frac{7\pm\sqrt{561}}{2},2^{5},-1^{17}\},$  \vskip 0.20 cm

(viiii). $\begin{bmatrix}
                          J-I_a & J   & 0 & 0 \\
                          J & R_{2k} &  J-S_{2k} &  0 \\
                          0 & J-S^{T}_{2k} & 0 & J \\
                          0 & 0 & J & 0  \\
                       \end{bmatrix}
                      $ where $(a,k)=(3,4)$ and $(5,3)$\\
   with spectra $\{1\pm3\sqrt{5},2^{4},-1^{10}\}$ and $\{2\pm\sqrt{43},2^{3},-1^{10}\},$  \vskip 0.20 cm

\end{thm}

From Theorem \ref{one}, we see that $\mathcal{F}$ contains four infinite families and twenty sporadic graphs. From the given
spectra it follows straightforwardly that

\begin{cor}\label{two} No two graphs  $ \mathcal{F}$ are cospectral.
\end{cor}

Given any two graphs $G$ and $H,$  let $G\cup H$ be the disjoint union of $G$ and $H,$ and $mG$  be the  disjoint union of $m$ copies of $G.$

\begin{thm}\label{four} Suppose $G$ and $G'$ are nonisomorphic cospectral graphs with at most two eigenvalues different from $2$, $-1$. Then $G=H \cup \beta K_{3}$ and $G'=H'\cup \beta' K_{3}$, where $H$ and $H'$ are
one of the following pairs of graphs in $\mathcal{F}:$

$(1). $ $H$ is of type (i) with $a=5$ and $k\geq2$, $H'$ is type (iv) with $k'\geq2$, where $5k=1+8k'$.

$(2). $ $H$ is of type (i) with $a=3$ and $k\geq2$, $H'$ is type (ii) with $k',\ell'\geq2$, where $k=k'\ell'$.

$(3). $ $H$ is of type (i) with $k\geq2$, $H'$ is type (viii) with $(a',k')=(4,10)$, where $a=1$ and $k=81$.

$(4). $ Both $H$ and $H'$ are of type (ii) with parameters $(k,\ell)$ and $(k',\ell')$, where $kl = k'\ell'$.
\end{thm}
\zm  The disjoint union of complete graphs in determined by its spectrum (see \cite{b2}). By Lemma \ref{nonsigular} (i), $G$ and $G'$ must have the described form.  Observing that $H$ and $H'$ has the eigenvalues $r > 2$ and $s < -1$, we easily find the given possibilities for $H$ and $H'$.
\ezm

It we take $\beta = 0$, we  can find the graphs in $\mathcal{F}$ having a non-isomorphic cospectral mate by Theorem \ref{four}. Hence, we have

\begin{cor}\label{fi} A graph $G\in\mathcal{F}$ is determined by its spectrum, unless $G$ is one of the following

$\diamond$ $G$ is of type (i) and $(a, k) = (1, 81)$.

$\diamond$ $G$ is of type (i) with $a = 3$ and $k$  is a composition number.

$\diamond$ $G$ is of type (i) with $a=5, k \equiv 5 \mod 8.$

$\diamond$ $G$ is of type (ii) and $k\ell$ has a divisor $d$ such that $\ell < d <k.$

\end{cor}

  By above Corollary \ref{fi}, then  the  generalized friendship graph $F_{t,r,k}$ with $t-r=3$
     is determined by its spectrum,  except when $r=1, k = 81; $ or $r=3,$  $k $ is  a composition number; or $r=5,
       k \equiv 5 \mod 8.$

\section{The proof of Theorem \ref{one}}
 In all cases  in Theorem  \ref{one}, we see that the corresponding quotient matrix has two eigenvalues different from $2$ and $-1$, and with Lemma \ref{odd} it straightforwardly follows that the remaining eigenvalues of the graph are all equal to $2$ and $-1$. So all graphs of Theorem \ref{one} are in ${\mathcal{F} }$.
\par  We choose $C$ to be a clique in $G\in \mathcal{F}$ with maximum size. By Lemma \ref{main} (graphs $G_1$ and $G_2$) $G$ contains no induced odd cycles of length five or more, therefore $|C| \geq 3.$ If there are more than one cliques of maximum size, we choose one for which the number of outgoing edges is minimal. The following lemmas and proposition are the key to our approach.

\begin{lem}\label{eight}
The vertex set of $C$ can be partitioned into two nonempty subsets $X$ and $Y$, such that the neighborhood of any vertex outside $C$ intersects $C$ in $X$, $Y$, or $\emptyset$.
\end{lem}
\zm The proof is analogous to the method in \cite{b5}. If $|C| = n - 1,$ the result is obvious. So assume $3 \leq |C| \leq n - 2$. Take vertices $x$ and $y$ outside $C$, and let $X$ and $Y$ consist of the neighbors of $x$ and $y$ in $C$, respectively.
Note that $X$ and $Y$ are proper subsets of $C$, since otherwise $C$ is not maximal. Suppose that $X \cap Y \neq\emptyset $, but $ X \nsubseteq Y $ . Then there exist vertices $ u \in X \cap Y $ and $ v \in X\backslash Y$. Let $w$ be a vertex in $ C \backslash X$. Then the subgraph induced by $\{u, v ,w , x , y\}$ is a forbidden subgraph $ G_3$, $G_4,$ or $G_5$. Therefore, if $X$ and $Y$ are not disjoint, then $ X \subseteq Y $, and analogously $ Y \subseteq X. $  Thus $ X \cap Y \neq \emptyset,$ implies $X = Y.$  If $ X \cap Y  = \emptyset $, assume there exist vertices $ u \in X $, $ v \in Y $, and $z \in C \backslash (X \cup Y )$, then $\{z, u, v, x, y\}$ induces a forbidden subgraph $G_6$ or $G_7$. This implies that if $X$ and $Y$ are disjoint and both nonempty, then $X \cup  Y = C$.
\ezm
\begin{lem}\label{nine}
If we take two vertices $x$ and $y$, $x\nsim y$,
 consider the corresponding  $2\times2$  principal submatrix $S$ of $A^{2}-A-2I$,
$$S=\left[
                        \begin{array}{cc}
                          d_{x}-2 & d_{xy}  \\
                          d_{xy}  &  d_{y}-2\\
                        \end{array}
                      \right].
 $$
 then $S$ is positive semi-definite and $\det S=(d_x-2)(d_y-2)-d_{xy}^{2}\geq0$.
\end{lem}
\par Let $\Gamma X$ and $ \Gamma Y$ denote the set of vertices outside $C$ adjacent to $X$ and $Y$ respectively. The set of vertices not adjacent to any vertex of $C$ will be denoted by $\Omega.$  Some of these sets may be empty, but clearly $\Gamma X$ or $\Gamma Y$ is nonempty (otherwise $G$ would be disconnected or complete). We choose $\Gamma X \neq\emptyset$  and distinguish three cases: (1) both $\Gamma Y$ and $\Omega$ are empty; (2) only $\Omega$ is empty; (3) $\Omega$ is nonempty. For convenience we define $a = |X|$, $b = |Y|$, and $c = |C| = a + b$. Let $G[ Z ] $ denote the induced subgraph by $ Z$.
\begin{prop} \label{there} Let $G$ be a graph, $ |X|=a$, $ |Y|=b$, $G[ \Gamma X] $ and  $G[ \Gamma Y] $ denote the induced subgraph by $ \Gamma X$ and $\Gamma Y$, respectively. Then

(i). If $b=1$ (resp., $  a=1$), then $G[\Gamma X]=lK_1$ (resp., $ G[\Gamma Y]=lK_1)$ ;

(ii). If $b=2$ (resp., $a=2)$, then $G[\Gamma X]=lK_1\cup kK_2$ (resp., $G[\Gamma Y]=lK_1\cup kK_2)$;

(iii).  If $b=3$ (resp., $a=3)$, then $G[\Gamma X]=lK_1\cup kK_2\cup mK_3$ (resp., $G[\Gamma Y]=lK_1\cup kK_2\cup mK_3)$;

(iv). If $b\geq 4$ (resp., $a\geq 4)$, then $G[\Gamma X]=lK_1\cup kK_2$ (resp., $G[\Gamma Y]=lK_1\cup kK_2)$.

\end{prop}
\zm (i). If $ b = 1$, then $ \Gamma X $ contains no edges, otherwise $ C $ would not be maximal.

(ii).   If $b=2$, choose  $ u \in X $, suppose $x \in \Gamma X$ has two neighbors $p$ and $q$ in $ \Gamma X$. If $p \nsim q$, then $\{ u, x, p, q,y\} $ ($y\in Y$) induces forbidden subgraph $G_3$ in Fig $2$, otherwise interchanging $\{x,p,q\}$ with $Y$ would give another larger clique. Therefore each vertex $x \in \Gamma X $ has at most one neighbor in $ \Gamma X $, and  $ G[\Gamma X]= lK_1\cup kK_2.$

(iii). If $b=3$, choose $ u \in X $, suppose $x \in \Gamma X$ has three neighbors $v $, $p$ and $q$ in $ \Gamma X$. If there exists a pair of vertex $p$ and $ q,$ such that  $p\nsim q,$ then $\{ u, x, p, q,y\} $  $(y\in Y)$ induces forbidden subgraph $G_3$, otherwise $v\sim p,$ $v\sim q,$ $p\sim q,$ interchanging $\{x,v,p,q\}$ with $Y$ would give another larger clique than before. Thus any vertex  of $ \Gamma X $ has at most two neighbor in $ \Gamma X $. If any vertex of $ \Gamma X$  has exactly two neighbor in $\Gamma X $, then the induced subgraph by $ \Gamma X$ are the disjoint union of cycles. If $G[\Gamma X]$ has a cycle with length four or more, then induces forbidden subgraph $G_3$, thus every cycle of length is three, and  $ G[\Gamma X]=lK_1\cup kK_2\cup mK_3.$

(iv).  If $b\geq 4$, let  $y, z, v ,w$ be four distinct vertices in $Y$, take a vertex $ u \in X $, suppose $x \in \Gamma X$ has two neighbors $p$ and $q$ in $ \Gamma X$. If $p\nsim q$, then $\{u, y, x, p, q \}$ induces forbidden  subgraph $G_3$, otherwise $\{u, y, z, v, w, x, p, q\}$ induces forbidden subgraph $G_8$. Thus each vertex $x \in \Gamma X $ has at most one neighbor in $ \Gamma X $,  and  $ G[\Gamma X]= lK_1\cup kK_2.$
   \ezm

\subsection {$\Gamma Y$ and $\Omega$ are empty}
\par  Assume that $1\leq b\leq3$, then $G[\Gamma X]=lK_1\cup kK_2\cup mK_3$ by Proposition \ref{there}. If $x \in \Gamma X $, $y\in Y$, then $d_{xy}=a$, $a\leq d_{x}\leq a+2$, $a\leq d_{y}\leq a+2$,  $\det S=(d_{x}-2)(d_{y}-2)-a^{2}\leq 0$. By Lemma \ref{nine}, $\det S=0$, thus $d_{x} = d_{y} = a+ 2$. Therefore $ G[\Gamma X]=mK_3$, $b=3$. Let $Y'=Y\cup \Gamma X=m'K_3$, $m'\geq 2$, since $Y$ and $\Gamma X$ are nonempty.  We can write $A$ as:
 $$ {\begin{matrix}
A=\begin{bmatrix}
J-I_a & J   \\
 J  & T_{3m'}\\
\end{bmatrix}
\end{matrix}}
$$
  where $3m'=|\Gamma X|+3,$ which leads to Case (i).

 Assume  that $b\geq 4$, then $G[\Gamma X]=lK_1\cup kK_2$  by Proposition \ref{there}.  By Lemma \ref{nonsigular} (ii), it is impossible that there exists one vertex of $\Gamma X$ has one neighbor in $\Gamma X$ but another vertex has no neighbor in $\Gamma X$.  We conclude that $G[ \Gamma X]=lK_{1} $ or $G[ \Gamma X]=kK_{2}$.

  Case (1):  $G[\Gamma X]=lK_{1}.$  If $l \geq 2$,  then there are at least two vertices have the same neighbors, which contradicts Lemma \ref{nonsigular} (ii). So $ l=1$  and we find
$$ {\begin{matrix}
A=\begin{bmatrix}\begin{smallmatrix}
J-I_a & J & J \\
 J & J-I_b  & 0 \\
  J & 0 & 0 \\
\end{smallmatrix}\end{bmatrix} ,&   Q = \begin{bmatrix}\begin{smallmatrix}
a-1 & b & 1  \\
 a & b-1 & 0\\
a &  0 & 0\\
\end{smallmatrix}\end{bmatrix}
\end{matrix}}.
$$
$P_Q(x)=a-ab-x+2ax+bx-2x^2+ax^2+bx^2-x^3$ shows that $Q$ has no eigenvalue $-1$ and has an  eigenvalue  $2$ if and only if $(a, b) = (7, 45),(8 ,27), (9, 21),(10,18),\\(12,15),(15,13),(18,12),(24,11)$ and $(42,10),$ which leads to Case (vi).
\par  Case (2): $G[ \Gamma X]=kK_{2}.$ If $k \geq 2$, then $G$ has eigenvalues 1, which contradicts Proposition \ref{p:wAPPROX}, thus $k=1$. $G$ has the following  $A$ and $Q$ :
$$ {\begin{matrix}
A=\begin{bmatrix}\begin{smallmatrix}
J-I_a & J  & J \\
 J & J-I_b  & 0 \\
  J & 0 & R_{2}\\
\end{smallmatrix}\end{bmatrix} ,&   Q = \begin{bmatrix}\begin{smallmatrix}
a-1 & b & 2  \\
 a& b-1 & 0 \\
a&  0 & 1 \\
\end{smallmatrix}\end{bmatrix}
\end{matrix}}.
$$
$P_Q(x)=1+a-b-2ab+x+2ax-x^2+ax^2+bx^2-x^3$ shows that $Q$ has no eigenvalue $-1$ and an  eigenvalue  $2$ if and only if $(a, b) = (2, 9),(3, 6)$ and  $(6, 5)$, which leads to Case (v).
\subsection{$\Gamma X$ and $\Gamma Y$ are nonempty, and $ \Omega $ is empty}
\subsubsection{Claim : $a \leq 3$ or $ b\leq 3$.}
\zm Suppose $ a \geq b \geq 4 $, by Proposition \ref{there}, we have $G[\Gamma X]=lK_1\cup kK_2$.   By Lemma \ref{nonsigular} (ii) and forbidden graphs $G_{20},G_{29},G_{30}$, it is impossible that there exists one vertex of $\Gamma X$ has one neighbor in $\Gamma X$ and another vertex has no neighbor in $\Gamma X$. We conclude that  $G[ \Gamma X]=kK_{2} $ or $G[ \Gamma X]=lK_{1}$. Forbidden graph $G_{28}$ implies that $k=1$. Similarly, we conclude that $G[ \Gamma Y]=K_{2} $, or $G[ \Gamma Y]=l'K_{1} $.
\par Case (1): $G[ \Gamma X]=K_{2}$, $G[ \Gamma Y]=K_{2} $.\\
Forbidden graph $ G_{20} $  implies that every vertex in $ \Gamma X $ is adjacent to  all vertices in $  \Gamma Y $. We
find the following $ A $ and $ Q $:
$$ {\begin{matrix}
A=\begin{bmatrix}\begin{smallmatrix}
J-I_a & J & J  & 0 \\
 J & J-I_b  & 0 & J \\
  J & 0  & R_{2} & J\\
  0 & J & J & R_2\\
\end{smallmatrix}\end{bmatrix} ,&   Q = \begin{bmatrix}\begin{smallmatrix}
a-1 & b & 2 & 0  \\
a& b-1 & 0 & 2\\
a &  0 & 1 & 2 \\
0 & b & 2 & 1\\
\end{smallmatrix}\end{bmatrix}
\end{matrix}}.
$$
 $P_Q(x)=-3+5a+5b-8ab-8x+5ax+5bx+4abx-6x^2-ax^2-bx^2-ax^3-bx^3+x^4$ shows that $Q$ has  no eigenvalue $-1$ and has  eigenvalue $2$ with multiplicity $1$ if and only if $(a, b) = (5, 4), $ but none of the other $ 3 $ eigenvalues are equal to $ 2 $ and $ -1 $. Thus the corresponding graphs are not in $\mathcal{F}$.

 Case (2): $G[ \Gamma X]=K_{2} $, $ G[\Gamma Y]=l'K_{1} $ .

Forbidden graph $ G_{29} $  implies that every vertex in $ \Gamma Y $ is adjacent to  all vertices in $  \Gamma X $. If  $l' \geq 2,$ then there are at least two vertices have the same neighbors, which contradicts Lemma \ref{nonsigular} (ii). So $l'= 1$, we find the following $ A $ and $ Q $:
$$ {\begin{matrix}
A=\begin{bmatrix}\begin{smallmatrix}
J-I_a & J & J  & 0 \\
 J & J-I_b  & 0 & J \\
  J & 0  & R_{2} & J\\
  0 & J & J & 0\\
\end{smallmatrix}\end{bmatrix} ,&   Q = \begin{bmatrix}\begin{smallmatrix}
a-1 & b & 2 & 0  \\
a& b-1 & 0 & 1\\
a &  0 & 1 & 1 \\
0 & b & 2 & 0\\
\end{smallmatrix}\end{bmatrix}
\end{matrix}}.
$$
$P_Q(x)=-2+2a+3b-3ab-5x+ax+3bx+3abx-3x^2-2ax^2-bx^2+x^3-ax^3-bx^3+x^4$ shows that $Q$ has no eigenvalue $-1,$ and has eigenvalue $2$ with multiplicity $1$ if and only if $(a, b) = (5, 5),$ but none of the other $ 3 $ eigenvalues are equal to $ 2 $ and $ -1 $. Thus the corresponding graphs are not in $ \mathcal{F} $.

 Case (3):  $ G[\Gamma X]=lK_{1} $,  $ G[\Gamma Y]=l'K_{1} $.\\
Now forbidden subgraph $ G_{30}$  implies that a vertex in $ \Gamma X $  is adjacent to all, or all but one vertices in $ \Gamma Y $, or all but two vertices in $  \Gamma Y $(and vice versa).
Let $ x $ be a vertex in $ \Gamma X $ and suppose $ x $ is adjacent to all vertices of $\Gamma Y $, suppose $ y $ is
another vertex in $ \Gamma X $, by Lemma \ref{nonsigular} (ii), $ y$ has fewer than $|\Gamma Y | - 4 $ neighbors in $ \Gamma Y $, contradiction. Similarly, if $ |\Gamma Y | \geq 2 $, then each vertex in $ \Gamma Y $ is adjacent to all but one vertices of $ \Gamma X $. This implies that the subgraph induced by $ \Gamma X \cup \Gamma Y $  is $ K_2 $ or a complete bipartite graph with the edges of a perfect matching deleted, by Lemma \ref{nonsigular} (ii), thus $l=l'$. Take two vertices $x'\in \Gamma X$, $y'\in \Gamma X,$ then $d_{x'}=d_{x'y'}+1$,  $d_{y'}=d_{x'y'}+1$, $\det S=(d_{x'}-2)(d_{y'}-2)-d_{x'y'}^{2}<0,$ by Lemma \ref{nine}, which is contradiction, therefore $l=l'\geq2$, the corresponding graphs are not in $ \mathcal{F} $.  We find  $G$ has the following $ A $ and $ Q $, where $l=l'=1$:
$$ {\begin{matrix}
A=\begin{bmatrix}\begin{smallmatrix}
J-I_a & J & J  & 0 \\
 J & J-I_b  & 0 & J \\
  J & 0  & 0 & 1\\
  0 & J & 1 & 0\\
\end{smallmatrix}\end{bmatrix} ,&   Q = \begin{bmatrix}\begin{smallmatrix}
a-1 & b & 1 & 0  \\
a& b-1 & 0 & 1\\
a &  0 & 0 & 1 \\
0 & b & 1 & 0\\
\end{smallmatrix}\end{bmatrix}
\end{matrix}}.
$$
$P_Q(x)=-1+a+b-ab-2x+2abx-2ax^2-2bx^2+2x^3-ax^3-bx^3+x^4$ shows that $Q$ has no eigenvalue $-1$, and  has eigenvalue $2$ with multiplicity $1$ if and only if $(a, b) = (9,9 ),(13,7),(21,6),$ but none of the other $ 3 $ eigenvalues are equal to $ 2 $ and $ -1 $. Thus the corresponding graphs are not in $ \mathcal{F} $.\ezm

\subsubsection{Claim  : $a \geq b=3.$}
\zm  First assume $a > b = 1$, by Proposition \ref{there}, we have $G[\Gamma X]=lK_1$. If $y \in  Y $  and $ x \in \Gamma X$, then $ x $  is adjacent to all vertices in $ \Gamma Y $, otherwise interchanging $ x $ and $ y $ would give another maximal clique of size $c$  with fewer outgoing edges. This implies that $x$  and $ y$ have the same neighbors, which is contradiction.
\par Next assume $a \geq b = 2 $, by Proposition \ref{there}, we have $G[\Gamma X]=lK_1\cup kK_2$.

 Suppose $G[\Gamma X]$ contains a $K_2$, then every vertex in $ \Gamma Y $ is adjacent to the two vertices of a $K_2$ in $ \Gamma X $. Otherwise interchanging two vertices of a $K_{2}$ in $G[\Gamma X]$ and $ Y $ would give another maximal clique of size $c$  with fewer outgoing edges. Choose a vertex $x$ of $K_2$ in $\Gamma X$, and a vertex $y$ of $Y$, thus $d_{x}=d_{y}=  d_{xy}+1$, $\det S=(d_{x}-2)(d_{y}-2)-d_{xy}^{2} <0$, which  contradicts Lemma \ref{nine}. Thus $G[\Gamma X]=lK_1.$

 Choose a isolated vertex $x$ of $\Gamma X$, for any vertex $y\in Y$, then  $d_{x}=d_{xy}$, by  Lemma  \ref{nonsigular} (ii), $d_{y}\geq d_{x}+5$ . If $a=2$, by Proposition \ref{there}, then $G[\Gamma Y]=k'K_2\cup l'K_1 $. By the same argument as above, we obtain $G[\Gamma X]=lK_1.$  Forbidden subgraph  $G_{12}$ shows that $d_{y}< d_{x}+5$, or we can find  two vertices $p,q\in \Gamma Y$, $d_p-d_q<5$,  which are contradiction. If $a\geq3$, then we have  $G[\Gamma Y]=l'K_1\cup k'K_2 \cup m'K_3$ by Proposition \ref{there}. Forbidden subgraphs  $G_{12},$  $G_{24}$, $G_{25}$, $G_{32}$ show that $d_{y}< d_{x}+5$, or we can find  two vertices $p,q\in \Gamma Y$, $p\nsim q$, $d_p-d_q<5$,  which are contradiction. \ezm

 We have $a \geq  b = 3 $, we have $G[\Gamma X]=mK_3\cup kK_2\cup lK_1$ by Proposition \ref{there}.

  Suppose $G[\Gamma X]$ contains a $K_2$, choose a vertex $x$ of $K_2$, for any vertex $y\in Y$, then  $d_{x}= d_{xy}+1$. By Lemma \ref{nine} $\det S=(d_x-2)(d_y-2)-d_{xy}^{2}\geq 0$, then   $d_{y}\geq  d_{xy}+4$. Forbidden subgraphs $G_{18}$, $G_{20}$  show that $d_{y}< d_{x}+4$, which is contradiction.

 Suppose $G[\Gamma X]$ contains a isolated vertex $x$, for any vertex $y\in Y$, then  $d_{x}=d_{xy}$, by  Lemma  \ref{nonsigular} (ii), $d_{y}\geq d_{x}+5$. But  forbidden subgraphs $G_{12}, G_{17}, G_{18},  G_{32}$ show that $d_{y}< d_{x}+5$, or we can find two vertices $p,q\in \Gamma Y$, $p\nsim q$, $d_p-d_q\leq 4$, contradiction.

  Thus $ G[\Gamma X]=mK_3,$  and every vertex in $ \Gamma X $ is adjacent to all vertices $ \Gamma Y $. Otherwise interchanging three vertices of a $K_{3}$ in $G[\Gamma X]$ and $ Y $ would give another maximal clique of size $c$ with fewer outgoing edges. By Lemma \ref{nonsigular} (ii), it is impossible that there exists one vertex of $\Gamma Y$ has no neighbor in $\Gamma Y$ but another vertex has one or two neighbor in $\Gamma Y$. If $G[\Gamma Y]= m'K_3\cup k'K_2$, then $a=3$, otherwise  $a\geq4$, which is impossible by forbidden subgraph $G_{8}$. Therefore  $a=3$, by the same argument as above, $ G[\Gamma Y ]=m'K_{3}$.  Thus $G[\Gamma Y ]=l'K_{1}$, $G[\Gamma Y ]=k'K_{2}$ or $ G[\Gamma Y ]=m'K_{3}$. Let $Y' = Y\cup \Gamma X=m''K_3$, then $m''\geq 2$, where $3m'' = |\Gamma X| + 3$, since $Y$ and $\Gamma X$ are nonempty.

  Case (1): $ G[\Gamma Y ]=l'K_{1}$ , if $l'\geq 2 $ , then there at least two vertices have the same neighbors, contradiction. So $l'= 1$ and  we find $G$ has the following $ A $ and $ Q $:
$$ {\begin{matrix}
A=\begin{bmatrix}\begin{smallmatrix}
J-I_a & J  & 0 \\
  J  & T_{3m''} & J\\
  0 & J &  0\\
\end{smallmatrix}\end{bmatrix} ,&   Q = \begin{bmatrix}\begin{smallmatrix}
a-1 &  3m'' & 0  \\
a& 2 & 1\\
0  & 3m'' & 0\\
\end{smallmatrix}\end{bmatrix}
\end{matrix}}.$$
 Computing det$(Q+I)$ and det$(Q-2I)$ shows that $Q$ has no  eigenvalues $-1$ and  $2$. Therefore the corresponding graphs are not in $ \mathcal{F} $.
\par Case (2):  $ G[\Gamma Y ]=k'K_{2}$,  $G$ has the following $ A $ and $ Q $:$$ {\begin{matrix}
A=\begin{bmatrix}\begin{smallmatrix}
J-I_a & J  & 0 \\
  J & T_{3m''} & J\\
  0  & J & R_{2k'}\\
\end{smallmatrix}\end{bmatrix} ,&   Q = \begin{bmatrix}\begin{smallmatrix}
a-1 &  3m'' & 0  \\
a& 2 & 2k'\\
0 &  3m'' & 1\\
\end{smallmatrix}\end{bmatrix}
\end{matrix}}.$$
 Computing $\det(Q+I)$ and $\det(Q-2I)$ shows that $Q$ has no  eigenvalues $-1$ and has an eigenvalue $2$  for $(a,k')=(6,1),(4,2)$, but $(a,k')=(4,2)$, $G$ has an eigenvalue 1, contradiction. Thus $(a,k')=(6,1),$ which leads to Case (iv).
\par Case (3):  $ G[\Gamma Y ]=m'K_{3}$, $a=3$. Let  $X'=X\cup \Gamma Y =lK_3$, then $l\geq 2$  as $X$ and $\Gamma Y$ are nonempty. Thus $G$ has the following $ A $ :
$$ A = \begin{bmatrix}
R_{3m''} & J   \\
J & R_{3l} \\
\end{bmatrix} $$
with $m'',l\geq2$, where $3m'' = |\Gamma X| + 3$ and $ 3l = | \Gamma Y | + 3$, which leads to Case (ii).

\subsection{ $ \Omega $ is nonempty }

Since $ G $ is connected, there exists an edge ${xz}$ with $z \in  \Omega $, and $x \in  \Gamma X$, or $x \in \Gamma Y.$ Assume
$ x \in \Gamma X$, take $u \in X$, and let $y$ be a neighbor of $z$ different from $x$. If $y \in  \Gamma Y$, then the
neighbor $v \in Y$  of $y$ together with $u,$ $x,$ $y,$ and $z$ induce a forbidden subgraph $G_1$ or $G_6$. Thus, $y \notin \Gamma Y$  which means $y\in \Gamma X \cup \Omega$. Similarly, if $x\in \Gamma Y$, then $y\in \Gamma Y\cup \Omega$. Without loss of generality, we assume that $\Gamma X$ and $\Omega$ are nonempty.

  \subsubsection{Claim : $a> b= 1$ or $a\geq b= 2$.}

\zm Assume $ a\geq b\geq 3$, it follows that $G[\Gamma X]= mK_3\cup kK_2\cup lK_1$  by Proposition \ref{there}. Forbidden subgraphs $G_{10},G_{19},G_{20}$ and Lemma \ref{nonsigular} imply that  at most one vertex in $\Omega$ is adjacent to all vertices in $ \Gamma X $. Similarly, at most one vertex in $\Omega$ is adjacent to all vertices in $ \Gamma Y $. Suppose $z\in \Omega$, then there is at least 2 vertices  in $\Gamma X$ by Lemma \ref{nonsigular} (i) , we can find two vertices $x$ and $y$, such that $x,y\in \Gamma X$, $x\sim z$, $y\sim z$.  Forbidden subgraphs $  G_{21}, G_{27}, G_{32}$  imply that every vertex in $ \Gamma X $ which is adjacent to an vertex of $\Omega$  has no neighbor in $ \Gamma X $, thus $x\nsim y$.  If $G[\Gamma Y]=\emptyset$, then $d_x=d_y=d_{xy}$; if $G[\Gamma Y]=l'K_1\cup k'K_2\cup m'K_3$, then  forbidden subgraph $G_{22}$ implies that   $x$ and $y$ are adjacent to all vertices of $K_2$ and $K_3$ in $ \Gamma Y $, forbidden subgraph $G_{11}$ implies that  isolate vertices in $\Gamma Y$ is adjacent to all vertices or all but one vertices in $ \Gamma X $ which is adjacent to $z$, forbidden subgraph $G_{13}$ implies that  a vertex in $\Gamma X$ which is adjacent to $z$ is adjacent to all vertices or all but one isolate vertex in $ \Gamma Y, $ thus $d_{xy}\leq d_x\leq d_{xy}+1$ and $d_{xy}\leq d_y\leq d_{xy}+1$, but $\det S\leq(d_{xy}-1)^{2}-d_{xy}^{2}<0$, which is contradicts Lemma \ref{nine}. Thus the corresponding graphs are not in $ \mathcal{F} $ for  $a\geq b\geq 3$. \ezm

  We have $a > b=1$ or $a \geq b =2$.
 \par If $a > b=1$, Then $G[\Gamma X] $ contains no edges, otherwise $ C $ is not maximal. Consider the set $Y' = Y\cup \Gamma X$, then $|Y'|\geq 2$, since $Y$ and $\Gamma X$ are nonempty. However $Y'$ contains no edges, otherwise  $C$ wouldn't be maximal. Let $Z$ be the set of vertices that are not in $X$ or $Y'$. Therefore $X$, $Y'$, and $Z$ give the following block structure of $A$:
 $$A=\left[
     \begin{array}{ccc} \begin{smallmatrix}
       J-I_a & J & 0 \\
       J & 0 & N \\
       0 & N^{T} & M \\
      \end{smallmatrix}\end{array}
   \right].$$

 Take three vertices $u \in X$, $x \in Y'$ and $y\in  Y'$. Consider the corresponding $3 \times 3$ principal submatrix $T$ of $A^{2}-A-2I$, then
$$ {\begin{matrix}
T=\begin{bmatrix} \begin{smallmatrix}
d_u-2&a-2&a-2&\\
a-2&d_x-2&d_{xy}\\
a-2&d_{xy}&d_y-2\\
 \end{smallmatrix}\end{bmatrix}
\end{matrix}}.
$$
Let  $T = (a-2)J + T'$, then
$$ {\begin{matrix}
T'=\begin{bmatrix}
d_u-a&0^{T}\\
0&T''\\
\end{bmatrix} ,&   T'' = \begin{bmatrix}
d_x-a&d_{xy}-a+2\\
d_{xy}-a+2&d_y-a\\
\end{bmatrix}
\end{matrix}}.
$$

Note that  $d_u> a$,  $d_x\geq a$ and  $d_y\geq a$. Without loss of generality, we  assume $d_y\geq d_x$. If $T''$ is positive definite, then so are $T'$ and $T$, which contradicts rank $T \leq2$. Therefore  $\det T''= (d_x-a)(d_y-a)-(d_{xy}-a+2)^{2} \leq0$, and by Lemma \ref{nine}  $\det S=(d_x-2)(d_y-2)-d^{2}_{xy}\geq 0$. If $d_{x} =  d_{xy}+1,$ then there exists $z$ such that $z\sim x$, but $z\nsim y$, then these neighbors of $y$ together with $x, y, z$ and any two vertices in $X$ induce forbidden subgraph $ G_{13}$, thus  $d_{xy}+1 \leq d_{y} \leq d_{xy}+ 3$, then $\det S\leq(d_{xy}-1)(d_{xy}+1)-d^{2}_{xy}< 0$, which is contradiction. If $d_{x} \geq  d_{xy}+2,$ then  det $T''> 0$, unless $d_{x} = d_{y} = d_{xy}+ 2$.  If $d_{x} =  d_{xy},$  then for any   two vertices $u, v $ of $Y'$ satisfy  $d_{u} = d_{v} = d_{uv}+ 2$ other than $x$. If  $|X|\geq 3$ and by Lemma \ref{nonsigular},  $d_y\geq d_x+5$,  which is impossible by forbidden subgraph $G_{31}$. If $|X|= a=2$, then det $T''= (d_x-a)(d_y-a)-(d_{xy}-a+2)^{2} =(d_x-2)(d_y-2)-d_{xy}^{2}\leq0$, and by Lemma \ref{nine}  $\det S=(d_x-2)(d_y-2)-d^{2}_{xy}\geq 0$, therefore det $T''$=det $S=(d_x-2)(d_y-2)-d^{2}_{xy}= 0$. Because $d_{x} =  d_{xy},$ then $d_{x} =3,$ $d_{y} =11,$ which is impossible by forbidden subgraph $G_9$; or $d_{x} =6,$ $d_{y} =11,$ which is impossible by forbidden subgraph $G_{16}$; or $d_{x} =4,$ $d_{y} =10,$ which is impossible by forbidden subgraph $G_{16}$, or $G$ has eigenvalue $1$, contradiction. Therefore, for any vertex of $Y'$, we conclude that $d_{x} = d_{y} = d_{xy}+ 2$, we find the following two possible structures for $N$:
$$N=\left[
       \begin{array}{ccc}
           J-S_{2k}^{T}   &0&J  \\
       \end{array}
     \right](k\geq2),~~ \text{or}~~       \\      N=\left[
       \begin{array}{ccc}
            S_{2m}^{T}  &0 &J\\
       \end{array}
     \right](m\geq3). $$
 Partition $Z=Z_1\cup Z_2\cup Z_3$ according to the structure on $N$, so that the vertices in $Z_2$ are not adjacent to all vertices of $Y'$ and $X$, the vertices in $Z_3$ are adjacent to all vertices of $Y'$. Forbidden subgraph $G_{13}$ implies that  $ G[Z_1]=mK_{2} $.  Suppose   $z\in Z_2$ is adjacent to a vertex of  $Z_1$, we can find  $u\in X$, $x, y\in Z_1$, $m,n\in Y'$, such that $m\sim x$,  $n\sim y$, $x\nsim y$, then these vertices $\{u, m, n, x, y, z\}$  induce forbidden subgraph $G_{2},$ or we find $u\in X$, $x, y\in Z_1$, $m\in Y'$, such that $m\sim x$,  $m\sim y$, $x\sim y$, $x\sim z$, then $\{u, m, x, y,  z\}$  induce forbidden subgraph $G_{7},$ thus the vertices in $Z_2$  are adjacent to all vertices of $Z_1$. Forbidden subgraph $G_{26}$ implies that at most one vertex in $Z_2.$ Suppose a vertex $z\in Z_3$ and $p\in Z_1$ are adjacent, we can find $u\in X$, $m,n\in Y'$, such that  $p\sim n$ and $p\nsim m$, then $\{u, m, n, z, p\}$ induce graph $G_{6}$, thus the vertices in $Z_3$  are non-adjacent to all vertices of $Z_1.$ Forbidden subgraph $G_{14}$ implies that every vertex in $Z_2$  is adjacent to all vertices of $Z_3.$ Forbidden subgraph $G_8$ implies that any vertex  of $ Z_3 $ has at most two neighbor in $Z_3 $. We can find two vertices $x'\in Z_1$, $y'\in Z_3$, $x'\nsim y'$, $d_{x'}=d_{x'y'}+1$, $d_{x'y'}+1\leq d_{y'}\leq d_{x'y'}+3$, det $ S=(d_{x'}-2)(d_{y'}-2)-d_{x'y'}^{2}<0$, by Lemma \ref{nine}, which is contradiction, therefore $Z_3$ is empty. Hence $N=[J-S_{2k}^{T}\   0]$ or $N=[S_{2k}^{T}\   0]$. Forbidden subgraph $G_{15}$ and Lemma \ref{nonsigular} imply that the second structures for $N$ is impossible. We find two structures for $Z_2$: $Z_2$ is empty, or $Z_2$ is nonempty and $|Z_2|=1$.

 Case (1): If $Z_2$ is empty,  $G[Y']=lK_{1}$, $ G[Z_1]=mK_{2} $, then $l=m$, and
 $G$ has the following adjacency matrix $A$ with quotient matrix $Q$:
 $$ {\begin{matrix}
A=\begin{bmatrix}\begin{smallmatrix}
J-I_a & J & 0 \\
 J & 0 & J-S_{2m}^{T}    \\
  0 &  J-S_{2m} & R_{2m} \\
\end{smallmatrix}\end{bmatrix} ,&   Q = \begin{bmatrix}\begin{smallmatrix}
a-1 &  m & 0 \\
a &  0 & 2m-2  \\
0& m-1 &1  \\
\end{smallmatrix}\end{bmatrix}
\end{matrix}}.
$$
 $P_Q(X)=2-2a-4m+3am+2m^2-2am^2+3x-ax-4mx+amx+2m^2x+ax^2-x^3$ shows that $Q$ has  no eigenvalue $-1$ and has an eigenvalue $2$   if and only if $(a, m) = (6, 3),(4 ,4)$, which lead Case (vii).

 Case (2): If $|Z_2|=1$, $G[Y']=lK_{1},$ $G[Z_1]=mK_{2},$ then $l=m$, and $G$ has the following adjacency matrix $A$ with quotient matrix $Q$:
 $$ {\begin{matrix}
A=\begin{bmatrix}\begin{smallmatrix}
J-I_a & J & 0& 0 \\
  J & 0 & J-S^{T}_{2m} & 0  \\
  0 & J-S_{2m} & R_{2m}  &J \\
  0& 0  &J & 0\\
\end{smallmatrix}\end{bmatrix} ,&   Q = \begin{bmatrix}\begin{smallmatrix}
a-1 & m & 0 & 0  \\
a &  0 &  2m-2  & 0 \\
0 &  m-1 &1 & 1 \\
0& 0& 2m& 0&  \\
\end{smallmatrix}\end{bmatrix}
\end{matrix}}.
$$
$P_Q(x)=(1+x)(2am^2-2x+2ax+2mx-amx-2m^2x-x^2-ax^2+x^3)$ shows that $Q$ has  an eigenvalue $-1$  and has an eigenvalue $2$   if and only if $a = 2$, we can rewrite $A$ as $$A=\left[
                   \begin{array}{cc}
                     R_{2m} & S-J_{2m} \\
                     J-S_{2m}^{T} & 0 \\
                   \end{array}
                 \right]
$$
with $m\geq 3$, which leads Case (iii).

  If $a\geq b=2,$ then $G[\Gamma X]= kK_2\cup lK_1$ by Proposition \ref{there}.  Forbidden subgraphs $G_{10},G_{19},G_{20}$ and Lemma \ref{nonsigular} imply that  at most one vertex in $\Omega$ is adjacent to all vertices in $ \Gamma X $. Similarly, at most one vertex in $\Omega$ is adjacent to all vertices in $ \Gamma Y $.
\par  Suppose $G[\Gamma X]$ contains a isolated vertex $x$, and $z\in \Omega$ is adjacent to $x$, choose any vertex $y\in Y$, then $ d_{x}= d_{xy}+1$. If $G[\Gamma Y]= \emptyset$, then $d_y=d_{xy}+1$; if $G[\Gamma Y]= l'K_1\cup k'K_2 \cup m'K_3$, then forbidden subgraph $G_{13}$ implies that  $d_{xy}+1\leq d_{y}\leq d_{xy}+3$, but det $S=(d_x-2)(d_y-2) -d^{2}_{xy}< 0,$ contradiction. Thus $z\nsim x$, choose a vertex $p\in \Gamma X$ of $K_2,$ such that $z\sim p$. Forbidden subgraph $G_{33}$ implies that $d_{x}= d_{xp}$. Therefore $d_x\geq 3$ by Lemma \ref{nonsigular}. But forbidden subgraphs $G_9, G_{16}$ imply that $d_p< d_x+5$, contradiction. Thus  the corresponding graphs are not in $ \mathcal{F} $ for   $G[\Gamma X]$ contains isolated vertices.

Thus $G[\Gamma X]=kK_2$. Consider the set $Y' = Y\cup \Gamma X=mK_2$, then $m\geq 2$, since $Y$ and $\Gamma X$ are nonempty. Let $Z$ be the set of vertices which are not in $X$ or $Y'$. Therefore $X$, $Y'$, and $Z$ give the following block structure of $A$:
 $$A=\left[
     \begin{array}{ccc}\begin{smallmatrix}
       J-I_a & J & 0 \\
       J & R_{2m} & N \\
       0 & N^{T} & M \\
     \end{smallmatrix}\end{array}
   \right].$$

   Take three vertices $u \in X$, $x \in Y'$ and $y\in  Y'$, $x\nsim y$. Consider the corresponding $3 \times 3$ principal submatrix $T$ of $A^{2}-A-2I$ , then
$$ {\begin{matrix}
T=\begin{bmatrix}\begin{smallmatrix}
d_u-2&a-1&a-1&\\
a-1&d_x-2&d_{xy}\\
a-1&d_{xy}&d_y-2\\
 \end{smallmatrix}\end{bmatrix}
\end{matrix}}.
$$
 Write $T = (a-1)J + T'$, then
$$ {\begin{matrix}
T'=\begin{bmatrix}
d_u-a-1&0^{T}\\
0&T''\\
\end{bmatrix} ,&   T'' = \begin{bmatrix}
d_x-a-1&d_{xy}-a+1\\
d_{xy}-a+1&d_y-a-1\\
\end{bmatrix}
\end{matrix}}.
$$
Note that $d_u>a+1$, $d_x\geq a+1$, $d_y\geq a+1$. Without loss of generality, we  assume $d_y\geq d_x$. If $T''$ is positive definite, then so are $T'$ and $T$, which contradicts rank $T \leq2$. Therefore det $T''= (d_x-a-1)(d_y-a-1)-(d_{xy}-a+1)^{2} \leq0$ and by Lemma \ref{nine} det $S=(d_x-2)(d_y-2)-d^{2}_{xy}\geq 0$. If $d_{x} =  d_{xy}+1,$  forbidden subgraphs $G_{20}, G_{23}$ show that  $d_{xy}+1 \leq d_{y} \leq d_{xy}+ 3$, then det $S=(d_x-2)(d_y-2)-d^{2}_{xy}< 0$, which is contradiction. If $d_{x} \geq  d_{xy}+2,$ then  det $T''> 0$, unless $d_{x} = d_{y} = d_{xy}+ 2$. We conclude that $d_{x} = d_{y} = d_{xy}+ 2$, we find the following two possible structures for $N$:
$$N=\left[
       \begin{array}{ccc}
           J-S_{2k}   &0 &J \\
       \end{array}
     \right](k\geq2), ~~\text{or}~~        N=\left[
       \begin{array}{ccc}
            S_{2m}
         &0&J\\
       \end{array}
     \right](m\geq3). $$
Partition $Z=Z_1\cup Z_2\cup Z_3$ according to the structure on $N$. Take five vertices $x,y\in Z_1$, $u\in X$,$m,n\in Y'$, such that $m\sim x$, $n\sim y$,  $m\nsim n$, if $x\sim y$ then $\{u,x,y,m,n\}$ induce graph $G_1$ in Fig $2$, thus  $ G[Z_1]=lK_{1} $. An argument  similar to the one used in $a>b=1$ shows that $Z_3$ is empty, and the second structures for $N$ is impossible. We find two structures for $Z_2$:  $Z_2$ is empty, or $Z_2$ is nonempty and $|Z_2|=1$.

 Case (1): If $Z_2$ is empty,   $G[Y']=kK_{2},$ $G[Z_1]=lK_{1},$ then $k=l$,  and
   $G$ has the following adjacency matrix $A$ with quotient matrix $Q$: $$ {\begin{matrix}
A=\begin{bmatrix}\begin{smallmatrix}
 J-I_a & J & 0 \\
J & R_{2k} &  J-S_{2k}  \\
 0 & J-S^{T}_{2k} & 0 \\
\end{smallmatrix}\end{bmatrix} ,&   Q = \begin{bmatrix}\begin{smallmatrix}
a-1 & 2k &  0 \\
a& 1 &  k-1  \\
0 &  2k-2 &0  \\
\end{smallmatrix}\end{bmatrix}
\end{matrix}}.$$
 $P_Q(x)=2-2a-4k+4ak+2k^2-2ak^2+3x-ax-4kx+2akx+2k^2x+ax^2-x^3$ shows that $  Q$ has  no eigenvalue $-1$ and  has  an eigenvalue $2$    if and only if $(a, k) = (4, 10),(5,7),(6,6),(9,5),$ which leads case (viii).
\par Case (2): If $|Z_2| =1$,   $G[Y']=kK_{2},$ $G[Z_1]=lK_{1},$ then $k=l$,  and
   $G$ has the following adjacency matrix $A$ with quotient matrix $Q$:
  $$ {\begin{matrix}
A=\begin{bmatrix}\begin{smallmatrix}
J-I_a & J & 0 &0 \\
  J   &R_{2k}& J-S_{2k}  &0   \\
  0 & J-S^{T}_{2k} & 0  &J\\
  0&0&J&0\\
\end{smallmatrix}\end{bmatrix} ,&   Q = \begin{bmatrix}\begin{smallmatrix}
a-1 &  2k & 0 & 0 \\
a  & 1& k-1 &0 \\
0 & 2k-2 &0 &1 \\
0&0&k&0\\
\end{smallmatrix}\end{bmatrix}
\end{matrix}}.
$$
$P_Q(x)=(1+x)(k-ak+2ak^2-2x+2ax+3kx-2akx-2k^2x-x^2-ax^2+x^3)$ shows that $Q$ has an eigenvalue $-1$  and has an eigenvalue $2$   if and only if $(a, k) = (3,4), (5,3), $ which leads case (viiii).\ezm

\vskip 0.4 true cm
\begin{center}{\textbf{Acknowledgments}}
\end{center}

This project was supported by the National Natural Science Foundation of China (No. 11571101).

 \end{document}